\documentclass[10pt]{article}

\title{Convex bodies 
in \\ Euclidean and  \weil geometries }
\author{Sumio Yamada}
\date{} 

 \usepackage{amsmath}
 \usepackage{amssymb}

\begin{document}
\newtheorem{theorem}{Theorem}
\newtheorem{proposition}[theorem]{Proposition}
\newtheorem{lemma}[theorem]{Lemma}
\newtheorem{corollary}[theorem]{Corollary}
\newtheorem{claim}[theorem]{Claim}
\newcommand{\weil}{Weil-Petersson }
\newcommand {\teich}{Teichm\"{u}ller }
\newcommand{\T}{{\mathcal T}}
\newcommand{\Tbar}{\overline{\mathcal T}}

\newenvironment{proof}{ {\sc {\bf Proof}}}{ {\sc {\bf q.e.d.}} \\}
\newenvironment{remark}{\noindent {\bf Remark}}{{\sc } \\}
\newenvironment{definition}{\noindent {\bf Definition}}{{\sc } }

\newcommand{\Tbdry}{\ensuremath{\partial {\cal T}}}
\newcommand{\g}{\ensuremath{\gamma}}
\newcommand{\del}{\partial}
\newcommand{\map}{\ensuremath{{\rm Map}(\Sigma)}}
\newcommand{\e}{\ensuremath{\varepsilon}}

\maketitle
\thanks{Supported in part by JSPS Grant-in-aid for Scientific Research No.20540201}

\section{Introduction}

Given an open convex body $\Omega$ in a Euclidean space, Hilbert in 1895 (\cite{H}) proposed a natural metric $H(x, y)$, now called Hilbert metric, defined on $\Omega$, as the logarithm of the cross ratio of a quadruple, $x, y, b(x, y)$ and $b(y,x)$, where $b(x, y)$ is where the ray from $x$ through $y$ hits the boundary of $\Omega$. 

For historical reasons, and in order to be consistent with the existing literature, we use the term ``metric"  in place of distance function; it does {\it not} refer to a Riemannian metric.  To be more precise, in this article we define a metric on a set $X$ to be
a function $\delta: X \times X \rightarrow ({\bf R}_+ \cup \{ \infty \} $ satisfying:
\begin{enumerate}
\item $\delta(x,x) = 0$ for all $x$ in $X$,
\item $\delta(x, z) \leq \delta(x, y) + \delta(y, z)$ for all $x,y$ and $z$ in $X$.
\end{enumerate}
In \cite{PT1}, for example, a function satisfiying the above is named weak metric. Note in this definition  neither the symmetry $\delta(x, y) = \delta(y,x)$ nor the nondegeneracy $\delta(x, y) = 0 \Rightarrow x=y$ is assumed.

  Then  the logarithm of the cross ratio as above indeed defines a 
metric, which is Finslerian and projective.  A Finsler structure on 
a Euclidean space determines a norm to each tangent space, and the 
norm itself is called Minkowski functional. A metric is said to be 
projective when Euclidean straight lines are geodesic.  The unit 
disc with its Hilbert metric $H(x, y)$ is a prominent example; 
Klein's model for the hyperbolic plane.

As the value $H(x, y)$ can be written as 
\[
\log \frac{|x-b(x,y)||y-b(y,x)|}{|y-b(x,y)||x-b(y,x)|} =  \log \frac{|x-b(x,y)|}{|y-b(x,y)|}
+ \log \frac{|y-b(y,x)|}{|x-b(y,x)|}
\]  
Funk~\cite{F} looked at the first term of the right hand side above as a metric, even though it is asymmetric, which has been called Funk metric.  The reader is referred to  \cite{PT1, PT2} where the historical and technical backgrounds are presented comprehensively.  

The purpose of this article is threefold.  The first is to introduce a new variational characterization ({\bf Theorem 1}) of the Funk metric for a convex domain $\Omega \subset {\bf R}^d$. It provides a new and short proof of the triangle inequality for the Funk and Hilbert distances. While the original proof by Hilbert uses the projective invariance of the cross ratio, the new one is less dependent on the Euclidean geometry, more on the convexity of the distance function. The new definition of the Funk metric also leads to new proofs for a set of results previously obtained in \cite{PT2} concerning the geometry of convex curves $\del \Omega$, as well as a new convexity statement of the functional $F(x_0, x)$ in the $x$ variable. Moreover we demostrate three different ways of writing down the Funk metric (and consequently the Hilbert metric) in this article. 

Secondly we generalize the Funk and Hilbert metrics  to \teich spaces of closed surfaces of genus $g >1$, where the background geometry is the one induced by the \weil metric (see, for example, \cite{W4}) instead of the Euclidean geometry.   The \weil geometry of \teich spaces has been known to possess many convex functionals (\cite{W1, Y0}) in addition to the \weil distance function itself. This feature has led to a construction of a new space (\cite{Y2}), called Teichm\"{u}ller-Coxeter complex, within which the original \teich space sits as an open convex subset.  This setting naturally allows us to adapt the Euclidean setting for the Funk metric ({\bf Theorem 2}) as well as the Hilbert metric to the \weil setting. In doing so, we introduce three different  (weak) metrics $F_1, F_2$ and $F_3$, each corresponding to one of the three representations of the Funk metric in the Euclidean setting, and each distinct from each other. The differences are caused by the negative sectional curvature of the \weil metric tensor. The \weil Funk metrics (as well as the \weil Hilbert metrics) thus defined are all invariant under the mapping class group action. 

The last aim is to draw the reader's attention to a similarity between the newly defined \weil Funk metric $F_2$ and so-called Thurston's asymmetric metric.  In his paper~\cite{Th}, Thurston embeds the \teich space to its cotangent space as a convex body, whose boundary is identified with the set of projective measured laminations on the given surfaces. There an asymmetric metric is introduced, by taking the logarithm of the ratio of hyperbolic lengths of a simple closed geodesic. In the last section of this article, we will demonstrate a similarity between the \weil Funk metric and Thurston's metric, which would lead to some questions and speculations.  When the closed surface is a torus $(g = 1)$, Thurston's metric has been analyzed closely in \cite{BPT} and it is effectively compared to the geometry of the \teich metric. We recall that on the \teich space for the torus, the \teich metric and the \weil metric are isometric.

The author would like to thank  Athanase Papadopoulos for many 
helpful comments, and the referee for the careful reading and 
the criticisms.     

\section{Euclidean Geometry}
\subsection{Funk metric}
Suppose that $\Omega$ is an open convex subset in a Euclidean space ${\bf R}^d$. In what follows, we set the presentation by Papadopoulos and Troyanov~\cite{PT1} as our reference for Funk 
and Hilbert metrics.

First we represent the convex set $\Omega$ as 
$
\cap_{\pi(b) \in {\cal P}} H_{\pi (b)}
$
where $H_{\pi(b)}$ is the half space bounded by a supporting
hyperplane $\pi(b)$ of $\Omega$ at the boundary point $b$, containing the convex set $\Omega$.
The index set $\cal P$ is the set of all supporting hyperplanes of $\Omega$. 
That for every boundary point $p$ there exists a supporting hyperplane $\pi(b)$ follows from the 
definition of convexity of $\Omega$.  In general, there can be more than one supporting hyperplane
of $\Omega$ at $p \in \del \Omega$. 
The index set $\cal P$ is identified with the set of unit normal vectors to the supporting hyperplanes.
It is identified with a subset of $S^{d-1}$, is equal to the entire sphere when the convex set is bounded.
We denote by ${\cal P}(b)$ the set of supporting hyperplanes at $b \in \del \Omega$.

\begin{definition}
For a pair of points $x$ and $y$ in $\Omega$, the Funk asymmetric metric~\cite{F} is defined by 
\[
F(x, y) =  \log \frac{d(x, b(x,y))}{d(y,  b(x,y))}. 
\]  
where the point $b(x,y)$ is the intersection of the boundary $\del \Omega$ and the ray $\{x+t \xi_{xy}\,\, | \,\, t>0\}$ 
from $x$ though $y$. Here $\xi_{xy}$ is the unit vector along the ray.
\end{definition}

Now let $\pi_0$ be a supporting hyperplane at $b(x,y)$, namely $\pi_0 \in {\cal P}(b(x,y))$.
Then  note the similarity of the triangle $\triangle (x, \Pi_{\pi_0}(x), b(x,y))$
and $\triangle (y, \Pi_{\pi_0}(y), b(x,y))$, where $\Pi_{\pi_0}(p)$ is the foot of the point $p$ on the hyperplane $\pi_0$, or put it differently $\Pi_{\pi_0} : {\bf R}^d \rightarrow \pi_0$ is the nearest point projection map. This says that 
\[
\log \frac{d(x, b(x,y))}{d(y,  b(x,y))} = \log \frac{d(x, \pi_0)}{d(y,  \pi_0)}.
\]
Also by the similarity 
argument of triangles, note that the right hand side of the equality is independent of the 
choice of $\pi_0$ in ${\cal P}(b(x, y))$.

Using the convexity of $\Omega$, the quantity $F(x, y)$ can be characterized variationally as follows. 
Define $T(x,\xi, \pi)$ by $\pi \cap \{x+t \xi  | t > 0\}$ with $\pi \in {\cal P}$. 
Consider the case $\xi = \xi_{xy}$.  When the hyperplane supports $\Omega$ at $p$, we have $T(x,\xi_{xy}, \pi) = b(x, y)$
and otherwise the point $T(x, \xi_{xy}, \pi)$ lies outside $\Omega$.  
When $\pi \notin {\cal P}(b(x,y))$,  by the similarity argument between the triangles
$\triangle (x, F_{\pi}(x), T(x, \xi_{xy}, \pi))$ and $\triangle (y, F_{\pi}(y), T(\xi_{xy}, \pi))$ again we have 
\[
\frac{d(x, \pi)}{d(y, \pi)} =\frac{d(x, T(x,\xi_{xy},\pi))}{d(y, T(x,\xi_{xy}, \pi))}.
\]
Note that the closest point to $x$ along the ray $\{x+t \xi_{xy} : t>0 \}$ of the form $T(x,\xi_{xy}, \pi)$ is $b(x, y)$.    
This in turn says that  a hyperplane $\pi$ which supports $\Omega$ at $b(x,y)$ maximizes the ratio $d(x, T(x,\xi_{xy}, 
\pi))/d(y, T(x,\xi_{xy}, \pi))$ among all the elements of ${\cal P}$;  

\[
 \log \frac{d(x, b(x,y))}{d(y,  b(x, y))} = \sup_{\pi \in {\cal P}} \log \frac{d(x, \pi)}{d(y,  \pi)}. 
\]

Hence we have a new characterization of the Funk metric;  
\begin{theorem} The Funk metric defined as above over a convex subset $\Omega \subset {\bf R}^d$ has the following variational formulation;
\[
F(x, y) = \sup_{\pi \in {\cal P}} \log \frac{d(x, \pi)}{d(y,  \pi)}.
\]    
\end{theorem}

With this formulation, one can readily see that $F(x, y)$ satisfies the triangle inequality, for
\[
F(x, y) + F(y, z)  =  \sup_{\pi \in {\cal P}} \log \frac{d(x, \pi)}{d(y, \pi)} +  \sup_{\pi \in {\cal P}} \log \frac{d(y, \pi)}{d(z, \pi)} 
\]
\[
  \geq   \sup_{\pi \in {\cal P}} \Big( \log \frac{d(x, \pi)}{d(y, \pi)} +   \log \frac{d(y, \pi)}{d(z, \pi)} \Big) =   \sup_{\pi \in {\cal P}} \log \frac{d(x, \pi)}{d(z, \pi)} = F(x, z)
\]

Note that the triangle inequality becomes an equality when 
\[
\sup_{\pi \in {\cal P}} \log \frac{d(x, \pi)}{d(y, \pi)} +  \sup_{\pi \in {\cal P}} \log \frac{d(y, \pi)}{d(z, \pi)} 
  =   \sup_{\pi \in {\cal P}} \Big( \log \frac{d(x, \pi)}{d(y, \pi)} +   \log \frac{d(y, \pi)}{d(z, \pi)} \Big)
\] 
is satisfied.  For this to occur, we only need ${\cal P}(b(x, y)) \cap {\cal P}(b(y,z)) \neq \emptyset$.  Let $\pi_0$ be an element in the set  ${\cal P}(b(x, y)) \cap {\cal P}(b(y,z)) \neq \emptyset$. Then 
the boundary point $b(x,y)$ and $b(y,z)$ share the same supporting hyperplane $\pi_0$, and we have  
\[
\sup_{\pi \in {\cal P}} \log \frac{d(x, \pi)}{d(y,  \pi)}
= \log \frac{d(x, \pi_0)}{d(y, \pi_0)}, \,\,\, 
\sup_{\pi \in {\cal P}} \log \frac{d(y, \pi)}{d(z,  \pi)}
= \log \frac{d(y, \pi_0)}{d(z, \pi_0)}
\]
and 
\[
\sup_{\pi \in {\cal P}} \log \frac{d(x, \pi)}{d(z,  \pi)}
= \log \frac{d(x, \pi_0)}{d(z, \pi_0)}.
\]
inducing the equality. The observation is equivalent to the following statement;
\\

\noindent {\bf Proposition 8.3} \cite{PT1} {\it Let $\Omega$ be an open convex subset of ${\bf R}^d$ such that $\del \Omega$ contains some non degenerate Euclidean segment $[p, q]$ and let $x$ and $z$ be two distinct points in $\Omega$ such that $\{t \xi_{xz}\,\, | \,\, t \geq 0\} \cap [p,q] \neq \emptyset$.  Let $\Omega'$ be the intersection of $\Omega$ with the affine subspace of ${\bf R}^d$ spanned by $\{x\} \cup [p,q]$.  Then, for any point $y$ in $\Omega'$ satisfying $\{t \xi_{xy}\,\, |t \,\, \geq 0\}  \cap [p,q] \neq \emptyset$ and $\{t \xi_{yz}\,\, | \,\, t \geq 0\}  \cap [p,q] \neq \emptyset$, we have $F(x, y) + F(y,z) = F(x,z)$.  }
\\

A notable situation when one has ${\cal P}(b(x, y)) \cap {\cal P}(b(y,z)) \neq \emptyset$ is when
$x,y$ and $z$ are colinear, with $y$ lying between $x$ and $z$.  This in turn says that
the straight lines are Funk geodesics, or in Hilbert's term, the Funk metric is projective.  This corresponds to Corollary 8.2 of \cite{PT1}.
Here a path $s:[a,b] \rightarrow (X, d)$ in a metric space $(X, d)$ is said to be 
geodesic when for any $a<t<b$, $F(s(a),s(t)) + F(s(t), s(b)) = F(s(a), s(b))$
is satisfied.

On the other hand, when $\pi_0$ is in the set  ${\cal P}(b(x, y)) \cap {\cal P}(b(y,z))$ the concatenation of the 
line segment $\overline{xy}$ and $\overline{yz}$ is also a Funk geodesic, a situation occurring when 
the boundary set $\del \Omega$ has a straight edge; a statement 
corresponding to Corollary 8.4 of \cite{PT1}.

We next consider the complementary situation where  ${\cal P}(b_1) \cap {\cal P}(b_2) = \emptyset$ for a 
pair of distinct points $b_i \in \del \Omega$. Geometrically this characterizes strict convexity of the domain $\Omega$, namely
the boundary $\del \Omega$ contains no closed line segments.  From the preceding argument, it follows that the only  
way the equality for the triangle inequality occurs is when the three points $x,y$ and $z$ are co-linear. Hence
for strictly convex domains, the Funk geodesics consists of line segments only, or equivalently, given a pair of points there is a unique Funk geodesic. This corresponds to Corollary 8.8 of \cite{PT1}.
    
With the additional assumption that the domain is strictly convex,
the following convexity of the Funk distance can be formulated.
Namely consider a geodesic/line $s(t)$ in $\Omega \subset {\bf R}^d$. Then the function $F(x, s(t))$ is convex in the parameter $t$;
first fix a supporting hyperplane $\pi$ in ${\cal P}$. Then it follows that
\[
\frac{d}{dt}  \log \frac{d(x, \pi)}{d(s(t), \pi)}=  -\frac{ \langle -\nu_{\pi} (s(t)), \dot{s}(t) \rangle}{ d(s(t), \pi ) } , \,\,\, 
\frac{d^2}{dt^2} \log \frac{d(x, \pi)}{d(s(t), \pi)} = \frac{\langle -\nu_\pi (s(t)), \dot{s}(t) \rangle^2}{d(s(t), \pi)^2}  
\]
Here $\nu_\pi(x)$ is the unit vector 
at $x$ perpendicular to the hypersurface $\pi(b)$ oriented toward $\pi$.  In particular $- \nu_b(x)$ is the gradient vector of the function $d(x, \pi)$

So $\displaystyle \log \frac{d(x, \pi)}{d(s(t), \pi)} $ is convex in $t$ for each $\pi \in {\cal P}$.  By taking supremum over $
\pi \in {\cal P}$, the resulting function $F(x, s(t))$ is convex in $t$. We note that the convexity of the geodesic ball $\{y \in \Omega \,\, | \,\, F_{\Omega} (x, y) < \delta \}$  centered at $x$  has been known (Proposition 8.11 (1) of \cite{PT1} ), which is equivalent to the convexity of $F(x, s(t))$ in $t$ for any lines $s(t)$ with $s(0) = x$.  

Next consider the function 
\[
F(s(t), x)= \sup_{\pi \in {\cal P}} \log \frac{d(s(t), \pi)}{d(x, \pi)}
\]
in $t$.  The first and second derivatives of this function are given by

\[
\frac{d}{dt}  \log \frac{d(s(t), \pi)}{d(x, \pi)}=  \frac{ \langle -\nu_\sigma (s(t)), \dot{s}(t) \rangle}{ d(s(t), \pi ) } , \,\,\, 
\frac{d^2}{dt^2} \log \frac{d(s(t), \pi)}{d(x, \pi)} = - \frac{\langle -\nu_\sigma(s(t)), \dot{s}(t) \rangle^2}{d(s(t), \pi)^2}  
\]

This says that the Funk distance $F(s(t), x)$ is not convex  in $t$.

\subsection{Finsler structure}

We next consider the linear structure of the Funk metric, by identifying it with a Finsler norm. 
The tautological weak Finsler structure~\cite{PT1} is given by the following Minkowski functional $p_{\Omega, x}(\xi) = \frac{1}{r_{\Omega, x}(\xi)}$ with
the radial function of $\Omega$ defined by 
\[
r_{\Omega, x}(\xi) = \sup \{t \in {\bf R} \,\, | \,\, x + t \xi \in \Omega \}.
\] 
From our exposition so far, we can identify the value of $r_{\Omega, x}(\xi)$ as 
\[
r_{\Omega, x}(\xi)= \inf_{\pi \in {\cal P}} d(x, T(x, \xi, \pi)) = \inf_{\pi \in {\cal P}} \frac{d(x, \pi)}{ \langle \nu_{\pi}(x), \xi \rangle}. \]

Note that the functional 
\[
p_{\Omega, x}(\xi) = \sup_{\pi \in {\cal P}} \frac{ \langle \nu_{\pi}(x), \xi \rangle}{d(x, \pi)}
\]
is convex in $\xi \in T_x {\bf R}^d$, as the functional $\displaystyle \frac{\langle \nu_{\pi}(x), \xi \rangle}{d(x, \pi)}$ is
convex  (linear in particular)  in $\xi$ for each fixed $\pi \in {\cal P}$, and by taking $\sup$ over $\pi $,  the 
convexity is preserved.   In Theorem 6.1 of \cite{PT1}, the infimum, among all the piecewise $C^1$ paths with given end points, of the length according to the Finsler norm  is shown to 
coincide with the Funk metric.   The proof follows from the fact that for a fixed $\pi$, the function
$\displaystyle \frac{ \langle \nu_{\pi}(\sigma(t)), \dot{\sigma}(t) \rangle}{d(\sigma(t), \pi)}$ is the $t$-derivative of
$\displaystyle \log \frac{d(x, \pi)}{d(\sigma(t), \pi)}$ for a $C^1$-path $\sigma(t)$.    

\subsection{A Brief Summary}
We have introduced three different representations of the Funk metric.  
The first one is the original definition by Funk~\cite{F};
\[
F_1(x, y) = \log \frac{d(x, b(x,y))}{d(y,  b(x,y))}
\]
where the point $b(x,y)$ is the intersection of the boundary $\del \Omega$ and the ray $\{x+t \xi_{xy}: t>0\}$ 
from $x$ though $y$.

The second one is the variational interpretation of the value above using the geometry of supporting hyperplanes;
\[
F_2 =  \sup_{\pi \in {\cal P}} \log \frac{d(x, \pi)}{d(y,  \pi)}
\]
where $\cal P$ is the set of all supporting hyperplanes of $\Omega$.   

And lastly  the Finsler structure $p_{\Omega, x}(\xi)$ is identified so that the Funk distance is described as the infimum of length of 
curves; 
\[
F_3 = \inf_{\sigma} \int_a^b p_{\Omega, \sigma (t)} (\dot{\sigma}(t)) dt 
\]
where the infimum is taken over all the piecewise $C^1$ curves with $\sigma(a) = x$ and $\sigma(b) = y$. 

We emphasize that they are all equal to each other 
\[
F(x, y) := F_1(x,y) = F_2(x,y)= F_3(x,y)
\]
for a convex domain $\Omega \subset {\bf R}^d$. 
 
\subsection{Hilbert metric}
We symmetrize the Funk metric by taking the arithmetic means
\[
H(x, y) = \frac12 \Big( F(x, y) + F(y,x) \Big)
\]
This metric is called Hilbert metric on $\Omega$.

For a convex set  $\Omega$ in ${\bf R}^d$, the affine line segment 
connecting $x$ and $y$ is a Funk geodesic realizing both lengths $F(x, y)$ 
and $F(y,x)$.  This says that the line segment also is the Hilbert geodesic. 

The Hilbert metric satisfies the triangle inequality which is obtained by taking the sum of the following pair of inequalities;
\[
F(x, y) + F(y,z) \geq F(x,z)
\]
which has been shown above and 
\[
F(y,x) + F(z, y) \geq F(z, x)
\]
which follows from the similar argument;
\[
F(y, x) + F(z, y)  =  \sup_{\pi \in {\cal P}} \log \frac{d(y, \pi)}{d(x, \pi)} + \sup_{\pi \in {\cal P}} \log \frac{d(z, \pi)}{d(y, \pi)} 
\]
\[
  \geq   \sup_{\pi \in {\cal P}} \Big( \log \frac{d(y, \pi)}{d(x, \pi)} +   \log \frac{d(z, \pi)}{d(y, \pi)} \Big)  = \sup_{\pi \in {\cal P}} \log \frac{d(z, \pi)}{d(x, \pi)} = F(z, x).
\]
Note this proof differs from Hilbert's which uses the projective property of the cross ratio among the
points $x,y, b(x,y)$ and $b(y,x)$.

Recall that for a strictly convex $\Omega$, $F(x, p(t))$ is convex in $t$ and the convexity of $F(p(t), x)$ is 
inconclusive in general.  As for the Hilbert metric where $F(x, p(t))$ and $F(p(t), x)$ are averaged, there 
is a complete characterization of $\Omega$ for which
it is convex. Suppose now that the Hilbert distance $H(x, p(t))$  of a given $\Omega$ is convex in $t$.  Then it 
is equivalent to saying that the geometry is non-positively  curved in the sense of Busemann, namely given a point $p$, there exists a 
neighborhood $U$ of $p$ such that for any pair $x$ and $y$, 
$
2 H(m_x, m_y) \leq H(x, y)
$
where $m_x$ is the midpoint with respect to $H$ between $x$ and $p$. Under this hypothesis
P.Kelly and E.G.Strauss~\cite{KS} showed that $\Omega$ is an ellipsoid, 
which is projectively equivalent to the unit disc. 
Hence the metric space $(\Omega, H)$ is isometric to the Poincar\'{e} disk.

\section{Weil-Petersson Geometry}

\subsection{$\Tbar$ as a convex subset in $D(\Tbar, \iota)$}

Now, in the \weil geometric setting, first recall that given a closed surface $\Sigma$ of genus $g >1$, 
the \weil metric completion(\cite{Y1}) $\Tbar$ is identified with
the augmented \teich space~\cite{Ab, M}, 
\[
\Tbar = \cup_{\sigma \in {\cal C}(\Sigma)} \T_\sigma  
\]
where ${\cal C}(\Sigma)$ is the curve complex of $\Sigma$, which can be identified as 
the index set of nodal surfaces $\Sigma_\sigma$'s specifying the location 
of the nodes. $\T_\sigma$ is the 
\teich space of the nodal surface $\Sigma_\sigma$. We also denote the set of all the elements $\sigma$  with $|\sigma|=1$ in ${\cal C}(\Sigma)$ by ${\cal S}$. Each element in $\cal S$ represent a single node.  In \cite{Y1}, the author has shown that $\Tbar$ is 
a CAT(0) space, and each frontier stratum $\Tbar_\sigma$ is embedded totally geodesically into $\Tbar$,
a generalization of a preceding result (\cite{W1}) by S. Wolpert which states that $\T$ is \weil geodesically convex.   
In \cite{Y2},  a new space has been introduced which can be viewed as a \weil geodesic completion,
called Teichm\"{u}ller-Coxeter complex $D(\Tbar, \iota)$. The space is a development of the original space $\Tbar$
by a Coxeter group generated by reflections across the frontier stratum $\{\Tbar_\sigma\}$.  It was shown by
Wolpert~\cite{W3} that two intersecting strata of the same dimension meet at a right angle (in the sense of 
the Alexandrov angle between \weil geodesics,) making the 
development $D(\Tbar, \iota)$ a so-called cubical complex~\cite{BH, D}. 
This feature then is used to show that the development too is a CAT(0) space. 
The Coxeter complex setting allows 
to view the \teich space as a convex set in an ambient space $D(\Tbar, \iota)$, bounded by a set of
 complex-codimension one ``supporting hyperplanes"  $\{D(\Tbar_\sigma, \iota) \,\, | \,\, |\sigma|=1\}$ of 
the frontier stratum $\{\Tbar_\sigma\}$ with each $\sigma$ 
representing a {\it single}  node. Every boundary point is contained in at least one of the set of the supporting hyperplanes $\{D(\Tbar_\sigma, \iota) \,\, | \,\, |\sigma|=1 \}$.  In this picture, each  $D(\Tbar_\sigma, \iota)$ is a totally geodesic
set, metrically and geodesically complete, and when $D(\Tbar_{\sigma_1}, \iota)$ and  $D(\Tbar_{\sigma_2}, 
\iota)$ intersect along  $D(\Tbar_{\sigma_1 \cup \sigma_2}, \iota)$, they meet at a right angle.  One can view this
as the translates of $\{D(\Tbar_\sigma, \iota)\,\, | \,\, |\sigma|=1\}$
by the action of the Coxeter group $W$ are forming a right-angled grid structure in $D(\Tbar, \iota)$,
whose lattice points are the orbit image by the Coxeter group $W$ of the set $\{\Tbar_\theta \,\, | \,\, |\theta|=3g-3 \}$ with $\theta$ indexing the maximal set of nodes on the surface.

Under this setting, for each $\sigma$ with $|\sigma| =1$ one can consider a {\it half space}, namely the set $H_{\sigma}$, containing $\Tbar$ in the
$D(\Tbar, \iota)$, and bounded by $D(\Tbar_\sigma, \iota)$. 
We note here the fact obtained by Wolpert~\cite{W3} that the \weil metric completion $\Tbar$ is the closure of the convex hull of the
vertex set  $\{\Tbar_\theta \,\, | \,\, |\theta|=3g-3 \}$, which suggests an interpretation of the \teich space 
as a simplex.  

We can summarize the discussion above as 
\[
\Tbar = \cap_{\sigma \in {\cal S}} H_{\sigma} \,\,\,  \mbox{ with } \,\,\, 
\partial \Tbar \subset \cup_\sigma D(\Tbar_\sigma, \iota).
\]
where every boundary point $b \in \partial \Tbar $ belongs to $D(\Tbar_\sigma, \iota)$ for some $\sigma$ in $\cal S$.  

We now transcribe the Euclidean Funk geometry as well as its compatible Finsler structure in the previous section to the \weil setting. 
There we demonstrated three different ways of writing down the Funk distance, which we called
$F_1,F_2$ and $F_3$, all of which coincide.  In the \weil setting, they differ from each other, and they are 
related by inequalities.

\subsection{\weil Funk metric $F_2$ }

First note that as each $\Tbar_\sigma$ lies in $\Tbar$ as a complete convex set, for each point $x \in \Tbar$, there 
exists the nearest point projection $\Pi_\sigma (x) \in D(\Tbar_\sigma, \iota)$, and the \weil geodesic $\overline{x 
\pi_\sigma(x)}$ meets with $D(\Tbar_\sigma, \iota)$ perpendicularly, its length uniquely realizing the distance
$\inf_{y \in \Tbar_\sigma} d(x, y) = d(x, \Pi_\sigma (x))$.  We denote this number by $d(x, \Tbar_\sigma)$. 
We also introduce a notation $\nu_\sigma(x)$   for the unit vector
at $x$ along the \weil geodesic  between $x$ and $\Pi_\sigma(x)$.  In particular $- \nu_\sigma(x)$ is the \weil 
gradient vector of the function $d(x, \Tbar_\sigma)$.  

Also note that one can replace $\Tbar_\sigma$ by $ 
D(\Tbar_\sigma, \iota)$ and get the same picture, for any $z \in \T$, we know that $\Pi_{\sigma}(z)$ is in $\T_{\sigma} \subset D(\iota, \Tbar_{\sigma})$
due to the fact that the frontier sets intersect perpendicularly in terms of the \weil Alexandorov angle. In other words, we have
$\Pi_\sigma^{-1} [\T_\sigma] = \T$ for any $\sigma \in {\cal C}(\Sigma)$.  
\\

\begin{definition}
We define the \weil Funk metric $F_2$ on $\T$ as;
\[
F_2(x, y) = \sup_{\sigma \in {\cal S}} \log \frac{d(x, \Tbar_\sigma)}{d(y, \Tbar_\sigma)}. 
\]
\end{definition}

In order to make the analogy with the Euclidean setting more obvious, and in order to make clearer the viewpoint that \teich space is a convex body within an ambient space, we can instead define the metric,  as
\[
F_2(x, y) = \sup_{\sigma \in {\cal S}} \log \frac{d(x, D(\Tbar_\sigma, \iota))}{d(y, D(\Tbar_\sigma, \iota))}
\]

The equality follows from the discussion in the paragraph preceding the definition of \weil Funk metric.

Note that $F_2(x, y)$ satisfies the triangle inequality, for
\begin{eqnarray*}
F_2(x, y) + F_2(y, z) & = & \sup_{\sigma \in {\cal S}} \log \frac{d(x, D(\Tbar_\sigma, \iota))}{d(y, D(\Tbar_\sigma, \iota))} +  \sup_{\sigma \in {\cal S}} \log \frac{d(y, D(\Tbar_\sigma, \iota))}{d(z, D(\Tbar_\sigma, \iota))} \\
 & \geq  & \sup_{\sigma \in {\cal S}} \Big( \log \frac{d(x, D(\Tbar_\sigma, \iota))}{d(y, D(\Tbar_\sigma, \iota))} +   \log \frac{d(y, D(\Tbar_\sigma, \iota))}{d(z, D(\Tbar_\sigma, \iota))} \Big) \\
   & = &  \sup_{\sigma \in {\cal S}} \log \frac{d(x, D(\Tbar_\sigma, \iota))}{d(z, D(\Tbar_\sigma, \iota))} = F_2(x, z)
\end{eqnarray*}

As in the Euclidean case, the equality holds when the stratum $D(\Tbar_\sigma, \iota)$ achieving the value of the 
supremum for $x$ and $y$ coincides with the one for $y$ and $z$.  However, those strata are no longer characterized 
geometrically, which prevents us from stating that the \weil geodesics are Funk distance minimizers. Recall in the 
Euclidean setting, the relevant supporting hypersurfaces are the ones containing the boundary points $b(x,y)$ and $b(y,z)$ 
where the rays hit the boundary set $\del \Omega$.   

We also remark that there is no definite statement for the convexity of the Funk distance $F_2(x, y)$ in the $y$ variable, unlike the Euclidean counterpart. 

We estimate the Funk distance $F_2(x, y)$ from above  and below by introducing two length metrics obtained by integrating  
linear structures $\tilde{p}_{\T, x}$ and  $\hat{p}_{\T, x}$ corresponding to Funk metrics $F_3$ and $F_1$ in our Euclidean setting.

\subsection{Finsler structure $\tilde{p}_{\T, x}$ and Funk metric $F_3$}

\begin{definition}
Define the first 
Finsler structure by the following Minkowski function:
\[
\tilde{p}_{\T, x} (\xi) =  \sup_{\sigma \in {\cal S}} \frac{ \langle \nu_{\sigma}(x), \xi \rangle_x}{d(x, \Tbar_\sigma)}
\]
where the inner product is the \weil pairing at $x$.
\end{definition}

Note that the function $\tilde{p}_{\T, x}$ is convex in $\xi \in T_x \T$, for it is the supremum of linear (in particular convex) functions. 
 Let $\alpha: [0,1] \rightarrow \T$ be a $C^1$-path connecting
$x = \alpha(0)$ and $y  = \alpha(1)$.  Then the length of $\alpha$ according to the Minkowski function $\tilde{p}_{\T, x}$ is defined as 
\[ 
\tilde{{\rm L}}(\alpha) =
\int_0^1 \tilde{p}_{\T, \alpha(t)} (\dot{\alpha}(t)) \,\, dt  
\]
This in turn defines a distance function
\[
F_3(x, y) =  \inf_{\alpha} \tilde{L}(\alpha)
\]
with $\alpha(0) = x, \alpha(1) = y$.

We show that the Funk distance $F_2(x, y)$ can be estimated from above by the length $\tilde{L}(\alpha)$;
\begin{eqnarray*}
F_2(x, y) & = & \sup_{\sigma \in {\cal S}} \log \frac{d(\alpha(0), \T_{\sigma})}{d(\alpha(1), \T_{\sigma})} \\
 & = & \sup_{\sigma \in {\cal S}} \Big[ - \log d(\alpha(t), \T_{\sigma}) \Big|_0^1 \Big] 
  =  \sup_{\sigma \in {\cal S}}  \int_0^1  \frac{ \langle \nu_{\sigma}(\alpha(t)), \dot{\alpha}(t) \rangle_{\alpha(t)}}{d(\alpha(t), \Tbar_\sigma)} \,\, dt  \\
  &\leq &   \int_0^1 \sup_{\sigma \in {\cal S}}   \frac{ \langle \nu_{\sigma}(\alpha(t)), \dot{\alpha}(t) \rangle_{\alpha(t)}}{d(\alpha(t), \Tbar_\sigma)}\,\, dt 
   =  \int_0^1 \tilde{p}_{\T, \alpha(t)} (\dot{\alpha}(t)) \,\, dt  = \tilde{{\rm L}}(\alpha)
\end{eqnarray*}

By taking infimum over the piecewise $C^1$ paths $\alpha$'s of the inequality, we obtain the estimate
\[
F_2(x, y) \leq F_3(x, y)
\]

\subsection{Linear structure $\hat{p}_{\T, x}$  and length metric $F_1$}

\begin{definition}
For a pair of distinct points $x$ and $y$ in $\T$, with the ray from $x$ through $y$ hitting a strata ${\T}_{\sigma_{xy}}$, we define  a number  $\phi_1(x, y)$ by 
\[
\phi_1(x, y)  =  \log \frac{d(x, T(x, \xi_{xy}, \sigma_{xy}))}{d(y, T(x, \xi_{xy}, \sigma_{xy}))}
\]
where 
$T(x, \xi, \sigma)$ is the point(s) where the ray $\{ \exp_x t \xi \,\,|\,\, t \geq 0 \}
$ and $D({\Tbar}_\sigma, \iota)$ meet and $\xi_{xy}$ is the tangent vector in $T_x \T$ such that $\exp_x \xi_{xy} = y$. If the ray $\{ \exp_x t \xi \,\,|\,\, t \geq 0 \}
$ does not hit any strata, then we set $\phi_1(x, y)=0$
\end{definition}  
\\

The ray staying in the interior $\T$ can indeed occur, when for example $x$ and $y$ are 
on a pseudo-Anosov axis~\cite{DW, W3} for some pseudo-Anosov 
element $\gamma$ in the mapping class group of $\Sigma$. 
It is not known that $\phi_1(x, y)$ satisfies the triangle 
inequality  satisfied by its Euclidean analogue $F_1(x, y)$.  \\

\begin{definition}
We now define another auxiliary linear structure on each tangent space;
\[
\hat{p}_{\T, x} (\xi) = \sup_{\sigma \in {\cal S}} \frac{\|\xi\|}{d(x, T(x, \xi, \sigma))}\]
where $\|\xi\|$ is the \weil norm of $\xi$. When the ray $t \xi$ hits no strata, set $\hat{p}_{\T, x} (\xi) = 0$
\end{definition}
\\

The functional $\hat{p}_{\T, x}$ is neither continuous nor convex 
on $T_x \T$, hence does not induce a Finsler structure.   
Once the geodesic ray in $D(\Tbar, \iota)$ crosses the bordification 
set $\partial 
\T \Tbar \backslash \T
$, the exponential map $\exp_x: T_x \T \rightarrow D(\Tbar, \iota)$ 
is no longer single valued. Hence for a given $\sigma$,  
$T(x, \xi, \sigma)$ is  possibly a set of points. However, in $
\Tbar$  the exponential map {\it is} single-valued.  Hence the 
value of 
$\sup_{\sigma \in {\cal S}} 1/d(x, T(x, \xi, \sigma))$ is the 
reciprocal of the \weil distance between $x$ and the bordification 
set $\partial 
\T $ along the ray $\{ \exp_x t \xi \,\,|\,\, t >0 \}$.  
Here we note when the ray from $x$ through $y$ hits the strata $\T_{\sigma_{xy}}$, the integral of the new linear structure along 
the \weil geodesic $\exp_x t \xi_{xy}$, which we denote $s(t)$ 
below, coincides with the number $\phi_1(x, y)$;
\[
\phi_1(x, y) =  \int_0^1 \hat{p}_{\Tbar, s(t)} (\dot{s}(t)) \,\, dt
\]
with $ s(1) = y$, as the $t$-derivative of $
\phi_1(x, s(t))$ is equal to 
$
 \hat{p}_{\T, s(t)} (\dot{s}(t)) 
$

This linear 
structure $\hat{p}_{\T, x}$ defines a length
for each piecewise $C^1$-path connecting
$x = \alpha(0)$ and $y  = \alpha(1)$;
\[
\hat{L}(\alpha) = \int_0^1 \hat{p}_{\T, \alpha(t)} (\dot{\alpha}(t)) \,\, dt  
\]
Then we have another weak metric $F_1(x, y)$ by
\[
F_1(x, y) = \inf_{\alpha} \hat{L}(\alpha) 
\]
where the infimum is taken over the set of piecewise $C^1$-paths connecting
$x = \alpha(0)$ and $y  = \alpha(1)$.  That this is a weak metric 
is clear as the two conditions 1 and 2 in the definition of metric 
are trivially met. When the geodesic ray from $x$ through $y$ stay in $\T$ indefinitely,
$F_1(x, y)=0$ as the ray has zero length with respect to the length defined by $\hat{p}_{\T, x}$.

We show that the Funk distance $F_2(x, y)$ can be estimated from below by $F_1(x, y)$.   As before, let $s(t) = \exp_x t \xi_{xy}$, and assume it hits $\T_{\sigma_{xy}}$, for otherwise $F_1(x, y)=0$ and the comparison is trivial.  
\begin{eqnarray*}
F_2(x, y) & = & \sup_{\sigma \in {\cal S}} \log \frac{d(\alpha(0), \T_{\sigma})}{d(\alpha(1), \T_{\sigma})}  \geq   \log \frac{d(\alpha(0), \T_{\sigma_{xy}})}{d(\alpha(1), \T_{\sigma_{xy}})} = \log \frac{d(x, \T_{\sigma_{xy}})}{d(y, \T_{\sigma_{xy}})}  \\
 & \geq & \log \frac{d(x, T(x, \xi, \sigma_{xy}))}{d(y, T(x, \xi, \sigma_{xy}))} = \phi_1(x, y) 
  =   \int_0^1  \frac{\| \dot{s}(t)  \|}{d( s(t), T(x, \xi, \sigma_{xy}))} \,\, dt \\
 & = & \int_0^1 \hat{p}_{\T, s(t)} ( \dot{s}(t) ) \,\, dt = \hat{L}(s)  
  \geq  \inf_{\alpha} \hat{L}(\alpha) \\
 &= & F_1(x, y) 
\end{eqnarray*}
The second inequality follows from the following argument.
First the \weil distance 
$d(\alpha(t), \T_\sigma)$ is strictly convex in $t$, which 
says that 
\[
\langle \nabla_{\dot{\alpha}(t)} \nu_\sigma(\alpha(t)), \dot{\alpha}(t) \rangle >0
\]
which implies that the angle between 
$\nu_\sigma$ and $\dot{\alpha}(t)$ is increasing in $t$, as the quantity
$\langle \nu_\sigma(\alpha(t)), \dot{\alpha}(t) \rangle$ is proportional to
the cosine of the angle.  Consider two right triangles the 
former being the geodesic triangle in $\Tbar$ whose vertices are $x$, $T(x, \xi, \sigma_{xy})$ and $\Pi_{\sigma_{xy}}(x)$ (the foot of $x$), and the latter being the Euclidean right triangle with a side of length $d(x, \Pi_{\sigma_{xy}}(x))$ between the vertex of a right angle and 
a vertex of an angle equal to the \weil angle at $x$ of the former
geodesic triangle in $\Tbar$.  

The monotone increasing property of the angle between 
$\nu_\sigma (\alpha(t)) $ and $\dot{\alpha}(t)$ contrasted to the 
constancy of the corresponding angles in the comparison triangle in ${\bf R}^2$, one concludes that 
\[
\frac{d(x, \T_{\sigma_{xy}})}{d(y, \T_{\sigma_{xy}})} > \frac{d(x, T(x, \xi, \sigma_{xy}))}{d(y, T(x, \xi, \sigma_{xy}))}.
\]
This completes the proof of the inequality $F_2 \geq F_1$.

Hence we summarize the \weil situation as follows.

\begin{theorem}
On the \teich space $\T$, the Funk metric $F_2$ and the Funk metrics $F_1$ and $F_3$ induced by $\hat{p}_{\T, x}$ and $\tilde{p}_{\T, x}$ respectively 
as above have the following comparison;
\[
F_1(x, y) \leq F_2(x, y) \leq F_3(x, y).
\]
Furthermore, each metric is invariant under the extended mapping class group action.
\end{theorem}

We recall that the extended mapping class group is the full \weil isometry group of $\T$(\cite{MW}.)

\subsection{\weil Hilbert metric on $\T$}
We symmetrize the \weil Funk metric $F_2$ by taking the arithmetic means.

\begin{definition}
We set 
\[
H(x, y) = \frac12 \Big( F_2(x, y) + F_2(y,x) \Big)
\]
and call this \weil Hilbert metric on the \teich space $\T$.
\end{definition}

This metric satisfies the triangle inequality.  It is obtained by taking the sum of the following pair of inequalities 
$
F_2(x, y) + F_2(y,z) \geq F_2(x,z)
$
which was shown above, and 
$
F_2(z, y)  + F_2(y, x) \geq F_2(z, x)
$
which follows from the same argument as in the proof of the triangle inequality 
in the Euclidean setting.

This provides a solution to a problem raised by A.Papadopoulos~\cite{P} {\it
``Realize Teichm\"{u}ller space as a bounded convex set somewhere and study the Hilbert metric on it" }
 as in our study, the \teich space is realized as a convex set, however unbounded, in the Teichm\"{u}ller-Coxeter space, on which we have defined a Hilbert metric.

We briefly remark that there are other ways of symmetrize the Funk metric, one of which is to define a new metric by
$D(x, y) = \max \{F(x, y), F(y, x) \}.$


\subsection{Concluding remarks}

On the closed surface $\Sigma$ as above, consider a pair of points in
the \teich space $\T$, represented by hyperbolic metrics $x$ and $y$.  
Then Thurston's asymmetric metric~\cite{Th} is defined by 
\[
K(x, y) = \sup_{\sigma \in {\cal S}} \log \frac{\ell_x (\sigma)}{\ell_y
(\sigma)}
\]
where $\ell_x(\sigma)$ denotes the hyperbolic length of the simple closed geodesic $\sigma$ on the hyperbolic surface $(\Sigma, x)$. The Finsler structure of this metric is described on p.20 of \cite{Th}; for  tangent vector $Y$ the Finsler norm is given by $\|Y\| =\sup_\lambda \Big( \frac{Y({\rm length} (\lambda)}{{\rm length}(\lambda)}  \Big)$ where $\lambda$ is taken over the space of measured geodesic laminations, corresponding to our \weil Finsler structure $\tilde{p}_{\T, x} (\xi)= \sup_{\sigma}  \frac{ \langle \nu_{\sigma}(x), \xi \rangle_x}{d(x, \Tbar_\sigma)}
$.   Note also that $\Big( \frac{Y({\rm length} (\lambda)}{{\rm length}(\lambda)}\Big)(x)$ is the directional derivative of $\log \frac{{\rm length}(\lambda)(y)}{{\rm length}(\lambda)(x)}$ by $Y$ at $x$, while $\frac{ \langle \nu_{\sigma}(x), \xi \rangle_x}{d(x, \Tbar_\sigma)}$ is the  directional derivative of $\log \frac{d(y, \Tbar_\sigma)}{d(x, \Tbar_\sigma)}$ by $\xi$ at $x$.

The asymptotic expansion (\cite{Y1}) of the \weil metric near 
a strata $\T_\sigma$  states that the \weil metric tensor  has an almost-product 
structure there, of the \weil metric of the frontier strata $\T_\sigma$ and the part expressed by the Fenchel-Nielsen parameters $(\ell_\sigma, \theta_\sigma)$.  The second part is of the form 
\[
4 \pi^3 \Big( d \rho_\sigma^2 + \frac14 \rho_\sigma^6 d \theta^2_\sigma \Big)
\] 
where $\rho_\sigma = \sqrt{\ell_\sigma}$. The constant $4 \pi^3$ is due to S.Wolpert(\cite{W3}.) This in turn says that near each stratum, the \weil distance $d(x, \T_
\sigma)$ is approximated by $\sqrt{\ell_\sigma(x)}$.  This observation can be transcribed to the following comparison;
\[
\log \frac{d(x, \Tbar_\sigma)}{d(y, \Tbar_\sigma)} \approx   \log \frac{\sqrt{\ell_x (\sigma)}}{\sqrt{\ell_y
(\sigma)}} = \frac12 \log \frac{\ell_x (\sigma)}{\ell_y
(\sigma)}
\]
provided that $x$ and $y$ are sufficiently close to $\Tbar_\sigma$ for some $\sigma$, an observation associating the two metrics $F_2$ and $K$.  

Most likely it is too optimistic to hope for the two metric $F_2$ and $K$ to coincide with each other (up to the constant $1/2$.)  However, this leads to a speculation as to 
whether the ratio $d(x, \Tbar_\sigma) / d(y, \Tbar_\sigma)$ can be seen as a number encoding some information about an appropriate geodesic lamination relevant to the pair of hyperbolic metrics $x$ and $y$ as well as the class $\sigma$ in $\cal S$.

Lastly, in addition to the similarity to Thurston's metric, the form 
of the \weil Funk metric $F_2$ should be contrasted with \teich 
metric defined on \teich spaces, in a particular representation 
found by Kerckhoff~\cite{Ke};
\[
d_T([G_1], [G_2]) = \frac12 \sup_{\sigma \in {\cal S}} \log \frac{{\rm Ext}_{[G_1]}(\sigma)}{{\rm Ext}_{[G_2]}(\sigma)}
\]
where ${\rm Ext}_{[G]}(\sigma)$ is the extremal length of $\sigma$ with respect to the conformal structure $[G]$.

\end{document}